\documentclass[12pt]{amsart}

\usepackage{amsfonts,amsmath,amsthm,amssymb,latexsym,amsthm,newlfont,enumerate,color,tikz-cd}

\newif\ifslide

\theoremstyle{plain}
\newtheorem{theorem}{Theorem}[section]

\newtheorem{lemma}[theorem]{Lemma}

\newtheorem{proposition}[theorem]{Proposition}
\newtheorem{definition-lemma}[theorem]{Definition-Lemma}

\newtheorem{red-question}[theorem]{\textcolor{red}{Question}}

\newtheorem{conjecture}[theorem]{Conjecture}

\theoremstyle{definition}
\newtheorem{definition}[theorem]{Definition}
\newtheorem{remark}[theorem]{Remark}

\newtheorem{example}[theorem]{Example}

\def\ideal#1.{I_{#1}}
\def\ring#1.{\mathcal {O}_{#1}}

\def\fring#1.{\hat{\mathcal {O}}_{#1}}
\def\proj#1.{\mathbb {P}(#1)}
\def\pr #1.{\mathbb {P}^{#1}}
\def\dpr #1.{\hat{\mathbb {P}}^{#1}}
\def\af #1.{\mathbb A^{#1}}
\def\Hz #1.{\mathbb F_{#1}}
\def\Hbz #1.{\overline{\mathbb F}_{#1}}
\def\fb#1.{\underset #1 {\times}}
\def\rest#1.{\underset {\ \ring #1.} \to \otimes}
\def\au#1.{\operatorname {Aut}\,(#1)}
\def\deg#1.{\operatorname {deg } (#1)}

\def\pic#1.{\operatorname {Pic}\,(#1)}
\def\pico#1.{\operatorname{Pic}^0(#1)}
\def\picg#1.{\operatorname {Pic}^G(#1)}
\def\ner#1.{NS (#1)}
\def\rdown#1.{\llcorner#1\lrcorner}
\def\rfdown#1.{\lfloor{#1}\rfloor}
\def\rup#1.{\ulcorner{#1}\urcorner}
\def\rcup#1.{\lceil{#1}\rceil}

\def\n1#1.{\operatorname {N_1}(#1)}  %Vector space of 1-cycles
\def\cn1#1.{\overline{\operatorname {N^1}(#1)}} %closure
\def\cone#1.{\operatorname {NE}(#1)}     %cone of curves
\def\ccone#1.{\overline{\operatorname {NE}}(#1)}
\def\none#1.{\operatorname {NF}(#1)}
\def\cnone#1.{\overline{\operatorname {NF}}(#1)}
\def\mone#1.{\operatorname {NM}(#1)} %cone of moving curves
\def\cmone#1.{\overline{\operatorname {NM}}(#1)}

\def\coef#1.{\frac{(#1-1)}{#1}}
\def\vit#1.{D_{\langle #1 \rangle}}
\def\mm#1.{\overline {M}_{0,#1}}
\def\H1#1.{H^1(#1,{\ring #1.})}
\def\ac#1.{\overline {\mathbb F}_{#1}}

\def\adj#1.{\frac {#1-1}{#1}}
\def\spn#1.{\overline{#1}}
\def\pek#1.#2.{\Cal P^{#1}(#2)}
\def\plk#1.#2.{\Cal P^{\leq #1}(#2)}
\def\ev#1.{\operatorname{ev_{#1}}}
\def\ilist#1.{{#1}_1,{#1}_2,\dots}
\def\bminv#1.{(\nu_1,s_1;\nu_2,s_2;\dots ;\nu_{#1},s_{#1};\nu_{r+1})}
\def\zinv#1.{(\nu_1,s_1;\nu_2,s_2;\dots ;\nu_{#1},s_{#1};0)}
\def\iinv#1.{(\nu_1,s_1;\nu_2,s_2;\dots ;\nu_{#1},s_{#1};\infty)}

\def\scr #1.{\mathcal #1}

\def\llist#1.#2.{{#1}_1,{#1}_2,\dots,{#1}_{#2}}
\def\ulist#1.#2.{{#1}^1,{#1}^2,\dots,{#1}^{#2}}
\def\lomitlist#1.#2.{{#1}_1,{#1}_2,\dots,\hat {{#1}_i}, \dots, {#1}_{#2}}
\def\lomitlistz#1.#2.{{#1}_0,{#1}_1,\dots,\hat {{#1}_i}, \dots, {#1}_{#2}}
\def\loc#1.#2.{\Cal O_{#1,#2}}
\def\fderiv#1.#2.{\frac {\partial #1}{\partial #2}}
\def\deriv#1.#2.{\frac {d #1}{d #2}}

\def\map#1.#2.{#1 \longrightarrow #2}
\def\rmap#1.#2.{#1 \dasharrow #2}
\def\emb#1.#2.{#1 \hookrightarrow #2}
\def\non#1.#2.{\text {Spec }#1[\epsilon]/(\epsilon)^{#2}}
\def\Hi#1.#2.{\text {Hilb}^{#1}(#2)}
\def\sym#1.#2.{\operatorname {Sym}^{#1}(#2)}
\def\Hb#1.#2.{\text {Hilb}_{#1}(#2)}
\def\Hm#1.#2.{\Hom_{#1}(#2)}
\def\prd#1.#2.{{#1}_1\cdot {#1}_2\cdots {#1}_{#2}}
\def\Bl #1.#2.{\operatorname {Bl}_{#1}#2}
\def\pl #1.#2.{#1^{\otimes #2}}
\def\mgn#1.#2.{\overline {M}_{#1,#2}}
\def\ialist#1.#2.{{#1}_1 #2 {#1}_2, #2\dots}
\def\pair#1.#2.{\langle #1, #2\rangle}
\def\vandermonde#1.#2.{\left|
\begin{matrix}
1 & 1 & 1 & \dots & 1\\
{#1}_1 & {#1}_2 & {#1}_3 & \dots & {#1}_{#2}\\
{#1}_1^2 & {#1}_2^2 & {#1}_3^2 & \dots & {#1}_{#2}^2\\
\vdots & \vdots & \vdots & \ddots & \vdots\\
{#1}_1^{#2-1} & {#1}_2^{#2-1} & {#1}_2^{#2-1} & \dots & {#1}_{#2}^{#2-1}\\
\end{matrix}
\right|
}
\def\vandermondet#1.#2.{\left|
\begin{matrix}
1 & {#1}_1   & {#1}_1^2 & \dots & {#1}_1^{#2-1}\\
1 & {#1}_2   & {#1}_2^2 & \dots & {#1}_2^{#2-1}\\
1 & {#1}_3   & {#1}_3^2 & \dots & {#1}_3^{#2-1}\\
\vdots & \vdots & \vdots & \ddots & \vdots\\
1 & {#1}_{#2}& {#1}_{#2}^2 & \dots & {#1}_{#2}^{#2-1}\\
\end{matrix}
\right|
}
\def\gr#1.#2.{\mathbb{G}(#1,#2)}

\def\alist#1.#2.#3.{{#1}_1 #2 {#1}_2 #2\dots #2 {#1}_{#3}}
\def\zlist#1.#2.#3.{#1_0 #2 #1_1 #2\dots #2 #1_{#3}}
\def\lomitlist30#1.#2.#3.{{#1}_0,{#1}_1 #2 \dots #2\hat {{#1}_i} #2\dots #2 {#1}_{#3}}
\def\lmap#1.#2.#3.{#1 \overset{#2}{\longrightarrow} #3}
\def\mes#1.#2.#3.{#1 \longrightarrow #2 \longrightarrow #3}
\def\ses#1.#2.#3.{0\longrightarrow #1 \longrightarrow #2 \longrightarrow #3 \longrightarrow 0}
\def\les#1.#2.#3.{0\longrightarrow #1 \longrightarrow #2 \longrightarrow #3}
\def\res#1.#2.#3.{#1 \longrightarrow #2 \longrightarrow #3\longrightarrow 0}
\def\Hi#1.#2.#3.{\text {Hilb}^{#1}_{#2}(#3)}
\def\ten#1.#2.#3.{#1\underset {#2}{\otimes} #3}
\def\lomitlist30#1.#2.#3.{{#1}_0 #2 {#1}_1 #2 \dots #2 \hat {{#1}_i} #2 \dots #2 {#1}_{#3}}
\def\mderiv#1.#2.#3.{\frac {d^{#3} #1}{d #2^{#3}}}

\def\Hom{\operatorname{Hom}}

\def\Supp{\operatorname{Supp}}

\def\Exc{\operatorname{Exc}}

\def\dim{\operatorname{dim}}

\def\deg{\operatorname{deg}}

\def\det{\operatorname{det}}

\def\Sing{\operatorname{Sing}}

\def\rest{\operatorname{res}}

\def\e{\Cal E}

\def\e1{E_1}
\def\e2{E_2}

\usepackage{hyperref,amscd,mathtools}

\DeclareOldFontCommand{\rm}{\normalfont\rmfamily}{\mathrm}

\author{Paolo Cascini}
\address{Department of Mathematics\\
Imperial College London\\
180 Queen's Gate\\
London SW7 2AZ, UK}
\email{p.cascini@imperial.ac.uk}

\author{Calum Spicer}
\address{Department of Mathematics, King's College London, Strand,
London WC2R 2LS, UK}
\email{calum.spicer@kcl.ac.uk}

\subjclass[2010]{14E30, 37F75}

\title{MMP for algebraically integrable foliations}
\dedicatory{Dedicated to Professor V.V. Shokurov on his 72nd birthday}

\begin{document}
\maketitle

\begin{abstract}
We show that termination of flips for $\mathbb Q$-factorial klt pairs in dimension $r$ 
implies existence of minimal models for algebraically integrable foliations 
of rank $r$ with log canonical singularities over a $\mathbb Q$-factorial klt projective variety. 
\end{abstract}
\tableofcontents

%!TeX root=main.tex

\section{Introduction}
Over the last few years, several results have appeared on the Minimal Model Program for foliations over a complex variety (e.g. see
\cite{mendes00,brunella00,McQ05,McQuillan08,spicer20,CS18, CS20}) and, in particular, it is known that foliations of rank one over any  manifold or foliations of any rank over a three-dimensional manifold  either are uniruled or they admit a minimal model. 

\medskip

The aim of this paper is to study the Minimal Model Program for algebraically integrable foliations with log canonical singularities. 
In our previous work, in collaboration with F. Ambro and V.V. Shokurov \cite{ACSS}, the cone theorem was proved for these foliations, and the MMP
was proven to exist for a special class, namely Property $(*)$ foliations.
While Property $(*)$ foliations are a very natural class of foliations, for many applications
Property $(*)$ is a too restrictive notion.  For instance, 
Fano foliations never satisfy Property $(*)$, cf. \cite[Theorem 5.1]{AD13}.

%\medskip

Our goal here is to show that, assuming termination of flips  for $\mathbb Q$-factorial klt pairs in dimension $r$,
the MMP exists for algebraically integrable 
foliations with log canonical singularities of rank $r$.  We remark that 
log canonical singularities form the largest class
of singularities where the MMP can be expected to be run. 

\medskip

The following is our main result: 
\begin{theorem}\label{t_main}
Assume termination of flips for $\mathbb Q$-factorial klt pairs in dimension $r$. 
Let $X$ be a normal $\mathbb Q$-factorial projective variety, 
let $\mathcal F$ be an algebraically integrable foliation on $X$ of rank $r$
and let $\Delta \geq 0$ be a $\mathbb Q$-divisor such that $(\mathcal F, \Delta)$ has log canonical singularities
%, $K_{\cal F}+\Delta$ is pseudo-effective
and $(X, \Delta)$ is klt. 

Then we may run a $(K_{\cal F}+\Delta)$-MMP $\varphi\colon X\dashrightarrow X'$ such that, if $\cal F'$ is the foliation induced on $X'$ and $\Delta'=\varphi_*\Delta$ then either
\begin{enumerate}
\item $(\cal F',\Delta')$ is a minimal model, i.e. $K_{\cal F'}+\Delta'$ is nef; or
\item $(\cal F',\Delta')$ admits a Mori fibre space, i.e. a contraction  $\eta\colon X'\to Z$ onto a normal projective variety $Z$ such that $\dim Z<\dim X$, $\rho(X/Z)=1$ and $-(K_{\cal F'}+\Delta')$ is relatively ample.  
\end{enumerate}
% $(\mathcal F, \Delta)$ admits a minimal model.
\end{theorem}

\medskip

\subsection{Sketch of the proof}
The main difficulty of the Theorem  is to prove the existence of the  contraction and, possibly, the flip associated to a $(K_{\cal F}+\Delta)$-negative extremal ray $R$, where   $(\cal F,\Delta)$ is an algebraically integrable log canonical foliated pair of rank $r$ over a $\mathbb Q$-factorial  variety $X$ such that $(X,\Delta)$ is klt (cf. Theorem \ref{thm_existence_contraction}). To this end, we divide the proof in three steps. We first run a suitable MMP to construct  a Property $(*)$-modification $\pi\colon \overline X\to X$ so that if  $\overline {\cal F}$ is the foliation induced on $\overline X$ and we write $K_{\overline{\cal F}}+\overline{\Delta}=\pi^*(K_{\cal F}+\Delta)$, then the pair $(\overline{\cal F},\overline {\Delta})$ satisfies Property $(*)$ (cf. Theorem \ref{t_starmod}). Then we run a partial $(K_{\overline{\cal F}}+\overline{\Delta})$-MMP so that we contract, and possibly flip, a $(K_{\overline{\cal F}}+\overline{\Delta})$-negative extremal ray $R'$ on $\overline X$ induced by $R$. Finally, we run a third MMP which contracts the strict transform of the divisorial part of the exceptional locus of $\pi$. We then show that the output of this third MMP is either the desired divisorial contraction or the flip associated to $R$. 
Note that we need termination of flips for $\mathbb Q$-factorial klt pairs in dimension $r$ only to run the second MMP. 

\subsection{Acknowledgements} We would like to thank F. Ambro and V.V. Shokurov for several discussions over some topics related to this paper.

%\begin{theorem}
%Assume termination of flips for $\mathbb Q$-factorial klt pairs in dimension $n-1$. 
%Let $X$ be a normal $\mathbb Q$-factorial projective variety of dimension $n$, 
%let $\mathcal F$ be an algebraically integrable foliation on $X$ of rank $\leq n-1$
%and let $\Delta \geq 0$ be a $\mathbb Q$-divisor so that 
%$(\mathcal F, \Delta)$ has log canonical singularities
%and $(X, \Delta)$ is klt. 
%
%Then $(\mathcal F, \Delta)$ admits a minimal model.
%\end{theorem}

%!TeX root=main.tex

\section{Preliminaries}

\subsection{Notations and preliminary results}
\label{s_preliminaries}
We work over an algebraically closed field of characteristic zero.
 
If $D$ is a $\mathbb Q$-divisor on a normal variety $X$ and $\Sigma$ is a prime divisor in $X$, then we denote by $m_\Sigma D$ the coefficient of $D$ along $\Sigma$. 
If  $D$ is a $\mathbb Q$-Cartier divisor on $X$, then we denote
\[
D^\perp\coloneqq \{\xi \in N_1(X)\mid D\cdot \xi=0\}.
\]
A morphism $f\colon X\rightarrow Z$ between normal varieties is said to be a {\bf contraction} if it is surjective and $f_*\mathcal O_X=\mathcal O_Z$.  
Given a contraction between normal varieties $f\colon X \rightarrow Z$ and an $\mathbb R$ divisor $B$
on $X$, we define its {\bf horizontal part} $B^h$ as the part of $B$ supported on prime divisors which dominate $Z$ and its {\bf vertical part} $B^v \coloneqq B-B^h$.

Given an equidimensional contraction $f\colon X \rightarrow Z$ between normal varieties and a $\mathbb Q$-divisor $B$ on $X$, we define the {\bf discriminant} of the pair $(X/Z,B)$ to be the $\mathbb Q$-divisor on $Z$ such that for any prime divisor $\Sigma$ on  $Z$ we have
\[m_\Sigma B_Z \coloneqq 1-\sup \{t \mid (\mathcal F, \Delta+tf^*\Sigma) \text{ is log canonical above the generic point of } \Sigma\}.\]

\medskip

We refer to \cite[Section 2.6]{CS18} for some of the standard definitions in the Minimal Model Program. Let 
$\varphi\colon X\dashrightarrow Y$ be a birational map between normal varieties.  The {\bf exceptional locus}  of $\varphi$ is the closed subset $\Exc\varphi$ of $X$ where $\varphi$ is not an isomorphism. We say that $\varphi$ is  a {\bf birational contraction} if $\Exc \varphi^{-1}$ does not contain any divisor. Let $\varphi\colon X\dashrightarrow Y$ be a birational contraction between normal varieties and let  $D$ be a $\mathbb Q$-Cartier $\mathbb Q$-divisor on $X$ such that $D_Y\coloneqq \varphi_*D$ is also $\mathbb Q$-Cartier. Then $\varphi$ is called {\bf $D$-negative} (resp. {\bf $D$-trivial}) if there exist a normal variety $W$ and birational morphisms $p\colon W\to X$ and $q\colon W\to Y$  resolving the indeterminancy locus of $\varphi$  such that $p^*D=q^*D_Y+E$ where  $E\ge 0$ is a $\mathbb Q$-divisor whose support coincides with the exceptional locus of $q$ 
(resp. $p^*D=q^*D_Y$). Note that by the negativity lemma (cf. \cite[Lemma 3.39]{KM92}), it follows that if $\varphi$ is a composition of steps of a $D$-MMP, then $\varphi$ is $D$-negative. 

\medskip

Given a divisor $G$ on a normal variety $X$, we denote by ${\rm Bs}(G)$ the base locus of the linear system $|G|$. If $D$ is a $\mathbb Q$-divisor on $X$, then the {\bf stable base locus}  of $D$ is
\[
\mathbb B(D)\coloneqq \bigcap {\rm Bs}(mD)
\]
where the intersection is taken over any sufficiently divisible positive integer $m$. The {\bf restricted base locus} of $D$ is
\[
\mathbb B_-(D)\coloneqq \bigcup_{\epsilon \in \mathbb Q_{>0}} \mathbb B(D+\epsilon A)
\]
where $A$ is a fixed ample divisor on $X$. Note that the definition does not depend on the choice of the ample divisor $A$. 

Similarly, we define the {\bf augmented stable base locus} of $D$ as 
\[
\mathbb B_+(D)\coloneqq \bigcap_{\epsilon \in \mathbb Q_{>0}} \mathbb B(D-\epsilon A)
\]
where $A$ is a fixed ample divisor on $X$. Also in this case, the definition does not depend on the choice of the ample divisor $A$.
By Nakamaye's theorem (e.g. see \cite[Theorem 1.4]{birkar17}), if $D$ is big and nef then the augmented locus of $D$ coincides with the {\bf exceptional locus} of $D$, which is 
 the union of all the subvarieties $V\subset X$ of positive dimension  such that $D|_V$ is not big.

\begin{lemma}\label{l_exceptional}
Let $\varphi\colon X\dashrightarrow Y$ be a birational contraction between normal $\mathbb Q$-factorial projective varieties and let $D$ be a  $\mathbb Q$-divisor on $X$ such that  $\varphi$ is $D$-negative.

Then $\Exc \varphi \subset \mathbb B_-(D)$. Moreover, if $\varphi_*D$ is nef, then the equality holds. 
\end{lemma} 

\begin{proof}
Let $D_Y\coloneqq \varphi_*D$. Since $\varphi$ is $D$-negative,  there exist a normal variety $W$ and birational morphisms $p\colon W\to X$ and $q\colon W\to Y$  resolving the indeterminancy locus of $\varphi$ and such that $p^*D=q^*D_Y+E$, where  $E\ge 0$ is a $\mathbb Q$-divisor whose support coincides with the exceptional locus of $q$. Since $\varphi$ is a birational contraction between $\mathbb Q$-factorial normal varieties, 
we have that $\Exc p\subset \Exc q$. It follows that $\Exc \varphi = p(\Exc q)$.

 Let $A$ be an ample divisor on $X$ and let $A_Y\coloneqq \varphi_*A$. 
 Then there exists $\epsilon>0$ sufficiently small, such that if we write 
 \[
 p^*(D+\epsilon A)=q^*(D_Y+\epsilon A_Y)+E',
 \]
 we have that $E'\ge 0$ and its support coincides with the support of $E$.  
 Thus, it follows that $\Exc q\subset \mathbb B(p^*(D+\epsilon A))$ and, therefore, $p(\Exc q)\subset  \mathbb B(D+\epsilon A)$. 
 Thus, $\Exc \varphi \subset \mathbb B_-(D)$, as claimed. 
 
 Let us assume now that $D_Y$ is nef. By the negativity lemma  (cf. \cite[Lemma 3.39]{KM92}), there exists an ample divisor $A_W$ on $W$ and a $\mathbb Q$-divisor $G\ge 0$, whose support is contained in $\Exc q$ and such that 
 $q^*A_Y=A_W+G$. Thus, $\mathbb B(q^*(D_Y+\epsilon A_Y))\subset \Exc q$. It follows that $\mathbb B(p^*(D+\epsilon A))\subset \Exc q$ and, therefore, $\mathbb B(D+\epsilon A)\subset \Exc \varphi$. Thus, our claim follows. 
\end{proof}

\begin{lemma}
\label{l_amplemodel}
Let $X\to U$ and $Y\to U$ be two projective morphisms of normal quasi-projective varieties. Let $D$ be a $\mathbb Q$-Cartier $\mathbb Q$-divisor on $X$ and let $\varphi\colon X\dashrightarrow Y$ be a $D$-trivial birational contraction over $U$  such that $\varphi_*D$ is ample over $U$. 

Then $\varphi\colon X\to Y$ is a morphism. 
\end{lemma}
\begin{proof}
Set $D_Y\coloneqq\varphi_*D$ and let $p\colon W\to X$ and $q\colon W\to Y$ be two birational morphisms from a normal variety $W$ which resolve the inderminancy locus of $\varphi$ and such that $p^*D=q^* D_Y$. Then $p^*D$ is semi-ample over $U$ and $H^0(X,\cal O_X(mD))\simeq H^0(Y,\cal O_Y(mD_Y))$ for any sufficiently divisible positive integer $m$.  It follows that also $D$ is semiample and that $\varphi$ is the morphism induced by $D$. 
\end{proof}

\medskip
\subsection{Foliations}
Let $X$ be a normal veriety. A {\bf foliation of rank $r$} on $X$ consists of a coherent subsheaf $T_\cal F \subset T_X$ of rank $r$ such that
\begin{enumerate}
\item $T_{\cal F}$ is saturated, i.e. $T_X/T_{\cal F}$ is torsion free; and

\item $T_{\cal F}$ is closed under Lie bracket.
\end{enumerate}
The {\bf co-rank} of $\cal F$  is its co-rank as a subsheaf of $T_X$.  We define the canonical divisor $K_{\cal F}$ to be a divisor on $X$ so that $\cal O_X(K_{\cal F}) \cong (\det T_{\cal F})^*$
and we define the conormal divisor $K_{[X/\cal F]}$ to be a divisor on $X$ so that
$\cal O_X(K_{[X/\cal F]}) \cong (\det (T_X/T_{\cal F}))^*$.
For any positive integer $d$, we denote $\Omega^{[d]}_X\coloneqq (\wedge^d\Omega^1_X)^{**}$.
Then, the foliation $\cal F$ induces a map 
\[
\phi\colon (\Omega^{[r]}_X \otimes \cal O_X(-K_{\cal F}))^{**} \rightarrow \cal O_X.\]
The {\bf singular locus} of $\cal F$, denoted by $\Sing \cal F$, is the cosupport of the image of $\phi$.

\medskip

Let $\sigma\colon Y\dashrightarrow X$ be a dominant map between normal varieties and let $\cal F$ be a foliation of rank $r$ on $X$. We denote by $\sigma^{-1}\cal F$ the {\bf induced foliation} on $Y$ (e.g. see \cite[Section 3.2]{Druel21}).  
If $T_{\cal F}=0$, i.e. if $\cal F$ is the foliation by points on $X$, 
then we refer to $\sigma^{-1}\cal F$ as the  {\bf foliation induced by $\sigma$} and we denote it by $T_{X/Y}$. 
In this case, the foliation $\sigma^{-1}\cal F$ is called {\bf algebraically integrable}. 
%In addition, if $\sigma$ is an almost holomorphic morphism, then we say that $\cal F$ is 
% {\bf non-dicritical}.

Let $f\colon X\to Z$ be a morphism between normal varieties and let $\cal F$ be the induced foliation on $X$. 
If $f$ is equidimensional,
%and $Z$ is smooth\footnote{P:do we need $Z$ smooth here? be don't assume it in Proposition \ref{prop_comparison} and I think that we don't need it if we assume that $f$ is equidimensional},
then we define the  {\bf ramification divisor} $R(f)$ of $f$  as 
\[
R(f)=\sum_D (\phi^*D-\phi^{-1}(D))
\]
where the sum runs through all the prime divisors of $Z$. In this case, we have
\[K_{\cal F} \sim K_{X/Z}-R(f)\] 
(e.g. see  \cite[Notation 2.7 and \S 2.9]{Druel15b}).

%Given a normal variety $X$ and 
%an algebraically integrable foliation $\mathcal F$ on $X$, we say that $\mathcal F$ is {\bf non-dicritical}
%if $\cal F$ is induced 
%for some foliated log resolution $p\colon X' \rightarrow X$ of $\mathcal F$,  
%such that $p^{-1}\mathcal F$ is induced by a contraction $X' \rightarrow Z$,  we have that the 
%induced rational map $X \dashrightarrow Z$ is almost holomorphic.
%\footnote{P: is this different from requiring that the foliation $\cal F$ is induced by an almost holomorphic rational map?}

Let $X$ be a normal variety and let $\cal F$ be a rank $r$ foliation on $X$. 
Let $S\subset X$ be a subvariety. Then  $S$ is said to be  {\bf $\cal F$-invariant}, or {\bf invariant by }$\cal F$, if for any open subset $U\subset X$ and any section $\partial \in H^0(U,\cal F)$, we have that 
\[ \partial (\mathcal I_{S\cap U})\subset \mathcal I_{S\cap U},
\]
where $\mathcal I_{S\cap U}$ denotes the ideal sheaf of $S\cap U$. 
If  $D\subset X$ is a prime divisor then we define $\epsilon(D) = 1$ if $D$ is not $\cal F$-invariant and $\epsilon(D) = 0$ if it is $\cal F$-invariant. Note that if $\cal F$ is the foliation induced by a dominant map $\sigma\colon Y\dashrightarrow Z$, then a divisor $D$ is $\cal F$-invariant if and only if it is vertical with respect to $\sigma$. 
%Given a divisor $D$ on $X$
%we may write $D = D_{{\rm inv}}+D_{{\rm n-inv}}$ where $D_{{\rm inv}}$ is the part of $D$ supported on $\mathcal F$-invariant
%divisors and $D_{{\rm n-inv}} \coloneqq D - D_{{\rm inv}}$.

\medskip

\subsection{Foliated pairs}
Let $X$ be a normal variety. 
A {\bf foliated pair} $(\cal F, \Delta)$ on $X$ consists of  a foliation $\cal F$ on $X$ and an $\mathbb R$-divisor 
  $\Delta\ge 0$
such that $K_{\cal F}+\Delta$ is  $\mathbb R$-Cartier.
Assuming that $\mathcal F$ is algebraically integrable, then we refer to \cite[Section 3.2]{ACSS} for some standard definitions of foliated singularities (e.g. foliated log smooth, log canonical, F-dlt). 

Let $(\mathcal F,\Delta)$ be a foliated pair on a normal variety $X$. Then, we say that $(\mathcal F,\Delta)$ satisfies {\bf Property $(*)$} if there exist a projective contraction $f\colon X\to Z$ and an $\mathcal F$-invariant divisor $G\ge 0$ on $X$ 
such that $\mathcal F$ is induced by $f$ and $(X/Z,\Delta+G)$ is a GLC pair which satisfies Property $(*)$ (cf. \cite[Section 2.2 and Definition 2.13]{ACSS}). In particular, there exists a  divisor $\Sigma_Z$ on $Z$ such that $(Z,\Sigma_Z$)  is log smooth and $G$ is the pre-image, as a reduced divisor, of $\Sigma_Z$ in $X$. 
 We say that such a divisor $G$ is  {\bf associated} to $(\mathcal F,\Delta)$ and we refer to $\Sigma_Z\subset Z$ as the {\bf discriminant} of $(\cal F,\Delta)$. By \cite[Lemma 2.14]{ACSS}, $\Sigma_Z$ coincides with the discriminant of the pair $(X/Z,\Delta+G)$ defined above. 
If $f\colon X\to Z$ is equidimensional then 
\[
K_{\cal F}+\Delta\sim_{f,\mathbb R} K_X+\Delta+G
\]
(cf. \cite[Proposition 3.6]{ACSS}).

\begin{proposition} \label{p_flipstar}
Let $X$ be a normal projective $\mathbb Q$-factorial klt  variety and 
let $(\mathcal F,\Delta)$ be a log canonical foliated pair on $X$ satisfying Property $(*)$. Assume that $\mathcal F$ is induced by an equidimensional projective contraction $f\colon X\to Z$. Let $R$ be a $(K_{\mathcal F}+\Delta)$-negative extremal ray. 

Then the contraction and, if the contraction is small, the flip associated to $R$ exist. Moreover, if $\varphi\colon X\dashrightarrow X'$  is the step of the MMP associated to $R$, $\mathcal F'$ is the induced foliation on $X'$ and $\Delta'=\varphi_*\Delta$, then $X'$ is $\mathbb Q$-factorial and  klt,  $(\mathcal F',\Delta')$ satisfies Property $(*)$ and $\mathcal F'$ is induced by an equidimensional projective contraction $f'\colon X'\to Z$. 
\end{proposition}
\begin{proof} We follow the same methods as in \cite[Lemma 3.14]{ACSS}. By the cone theorem (cf. \cite[Theorem 3.9]{ACSS}), we may assume that $R$ is spanned by a curve $\xi$ on $X$ which is tangent to $\mathcal F$.  Let $G\ge 0$ be a divisor associated to $(\mathcal F,\Delta)$. Then \cite[Proposition 3.6]{ACSS} implies
\[
K_{\cal F}+\Delta\sim_{f,\mathbb R} K_X+\Delta+G,
\]
and, by \cite[Lemma 2.14]{ACSS}, the pair $(X,\Delta+G)$ is log canonical. It follows that $R$ is $(K_X+(1-\epsilon)\Delta+(1-\epsilon)G)$-negative for some sufficiently small $\epsilon>0$ such that  $(X,(1-\epsilon)\Delta+(1-\epsilon)G)$ is klt. Thus,
 the existence of the contraction follows from the classic base point free theorem  (cf. \cite[Theorem 3.3]{KM98}). Similarly, if the contraction is small, then the flip associated to $R$ exists by the classical result on existence of flips. In particular, $X'$ is $\mathbb Q$-factorial and  klt. 

Finally, the last part of the Proposition follows by \cite[Proposition 2.18]{ACSS}.
\end{proof}

The following is a slightly stronger version of \cite[Theorem 3.10]{ACSS}:
\begin{theorem}\label{t_starmod}
Let $(\mathcal F,\Delta)$ be a log canonical foliated pair on a normal projective pair $X$ and assume that  $\mathcal F$ is algebraically integrable. 

Then there exists a birational morphism $\pi\colon X'\to X$ such that if $\mathcal F'$ is the foliation induced on $X'$ and $\Delta'\coloneqq \pi_*^{-1}\Delta+\sum \epsilon(E)E$, where the sum runs over all the exceptional divisors of $\pi$, then 
\begin{enumerate}
\item $X'$ is klt and $\mathbb Q$-factorial;
\item $\mathcal F'$ is induced by an equidimensional morphism $f\colon X'\to Z$; 
\item $(\mathcal F',\Delta')$ satisfies Property $(*)$ and if $G$ is the associated divisor, then every $\cal F'$-invariant divisor in the exceptional locus of $\pi$ is contained in the support of $G$; and 
\item $K_{\mathcal F'}+\Delta'=\pi^*(K_{\cal F}+\Delta)$.
%\item every divisor in the exceptional locus of  $\pi$ is $\mathcal F'$-invariant. 
\end{enumerate}
\end{theorem}
The morphism $\pi$ is called a {\bf Property $(*)$-modification}. 
\begin{proof} 
By \cite[Theorem 2.2 and Proposition 2.16]{ACSS}, 
there exists a birational morphism  $m\colon W\to X$  such that, if  $\cal F_W$ is the foliation induced on $W$ and 
\[
\Delta_W\coloneqq m^{-1}E+\sum \epsilon(E)E
\]
where the sum runs over all the $m$-exceptional divisors, then $\cal F_W$ is induced by an equidimensional contraction  $g\colon W\to Z'$
and $(\cal F_W,\Delta_W)$ satisfies Property $(*)$. Let $G_W$ be its associated divisor. We may assume that the support of $G_W$ contains every $\cal F_W$-invariant $m$-exceptional divisor.  By Proposition \ref{p_flipstar}, we may run a $(K_{\mathcal F_W}+\Delta_W)$-MMP over $X$. By \cite[Lemma 3.13]{ACSS}, every curve contracted by this MMP is also mapped to a point in $Z'$ and, therefore, every step of the MMP is also a step in a relative $(K_W+\Delta_W+G_W)$-MMP over $X$. In particular, the MMP terminates. Call  this MMP $\varphi\colon W\dashrightarrow X'$  and let $\pi\colon 
X'\to X$ be the induced birational morphism. Let $\cal F'$ be the foliation induced on $X'$, $\Delta'\coloneqq \varphi_*
\Delta_W$  and $G\coloneqq \varphi_*G_W$. Proposition \ref{p_flipstar} implies that   $(\cal F',\Delta')$ satisfies Property $(*)$ and if $G$ is its associated divisor, 
then the support of $G$ contains every $\mathcal F'$-invariant divisor in the exceptional locus of $\pi$. Thus, (3) follows.  In addition, Proposition \ref{p_flipstar} implies (1) and (2). Finally, by construction, $\varphi$ contracts every $m$-exceptional  divisor $E$ such that $a(E,\mathcal F,\Delta)>-\epsilon (E)$. Thus, (4) follows. 
\end{proof}

\begin{lemma}\label{l_termination*}
Assume termination of flips for $\mathbb Q$-factorial klt pairs in dimension $r$. Let $X$ be a normal projective $\mathbb Q$-factorial klt variety and let $(\cal F, \Delta)$ be a log canonical foliated pair of rank $r$ on $X$ satisfying Property $(*)$. Assume that $\cal F$ is induced by an equidimensional projective contraction $f\colon X\to Z$. 

Then any sequence of $(K_{\cal F}+\Delta)$-flips terminates. 
\end{lemma}

\begin{proof}
Assume by contradiction that  
\[
X_1\coloneqq X\dashrightarrow X_2\dashrightarrow \dots
\]
is an infinite sequence of $(K_{\cal F}+\Delta)$-flips. 
By Proposition \ref{p_flipstar}, it follows that each $X_i$ is $\mathbb Q$-factorial and klt. Moreover, if $\cal F_i$ denotes the foliation induced on $X_i$ and  $\Delta_i$ is the strict transform of $\Delta$ on $X_i$, then $(\cal F_i, 
\Delta_i)$ satisfies Property $(*)$. Let $G_i$ be the associated  divisor to $(\cal F_i,\Delta_i)$. 
As in the proof of Proposition \ref{p_flipstar}, it follows that  
each map $\varphi_i\colon X_i\dashrightarrow X_{i+1}$  is a $(K_{X_i}+\Delta_i+G_i)$-flip over $Z$ 
and the induced morphism $f_i\colon X_i\to Z$ is equidimensional. By Noetherian 
induction, there exists a closed point $p\in Z$ such that, if we denote by $F_i$ the fibre $f_i^{-1}(p)$, then the induced map $F_i\dashrightarrow F_{i+1}$ is not an isomorphism for infinitely many $i$. Note that, by assumption, $F_i$ has dimension $r$.  By \cite[Proposition 2.18]{ACSS}, it follows that $\varphi_i^{-1}$ is an isomorphism at the generic point of each irreducible component of $F_{i+1}$. Thus, we may assume that there exists an irreducible component $S_i$ such that the restriction $\psi_i\coloneqq \varphi_i|_{G_i}\colon S_i\dashrightarrow S_{i+1}$ is a birational map which is not an isomorphism for any sufficiently large integer $i$. 
Let $\Sigma_Z$ be the boundary 
associated to $(\cal F, \Delta)$ and let $H$ be a sufficiently general divisor on $Z$ such that $(Z,\Sigma_Z+H)$ is log smooth and $p$ is a zero dimensional stratum of $\Sigma_Z+H$. By the definition of Property $(*)$, it follows that 
$(X_i,\Delta_i+G_i+f_i^*H)$ is log canonical around $S_i$. 
Note that each $X_i\dashrightarrow X_{i+1}$ is a $(X_i,\Delta_i+G_i+f_i^*H)$-flip 
and that \cite[Lemma 2.15]{ACSS} implies that $S_i$ 
is a log canonical centre for $(X_i,\Delta_i+G_i+f_i^*H)$. 
Let $S_i^\nu$ the normalisation of $S_i$. 
We have that $F_i$ can be written as the intersection of $\dim X-r$ 
log canonical centres of $(X_i,\Delta_i+G_i+f_i^*H)$ of codimension one. 
Thus, inductively applying adjunction, 
  we may write 
\[
(K_{X_i}+\Delta_i+G_i+f_i^*H)|_{S_i^{\nu}}=K_{S_i^{\nu}}+\Delta_{S_i}
\]
for some $\mathbb R$-divisor $\Delta_{S_i}\ge 0$ on $S_i^{\nu}$. 
By proceeding as in the proof of Special Termination (e.g. see \cite[Theorem 4.2.1]{Fujino05}), 
we can show that there exists a  sufficiently large $i_0$ and a $\mathbb Q$-factorial 
dlt modification $\mu\colon T\to S_{i_0}^{\nu}$  so that, if we write 
$\mu^*(K_{S_{i_0}}^{\nu}+\Delta_{S_{i_0}})=K_{T}+\Delta_T$,  
then the maps 
\[
S_{i_0}\dashrightarrow S_{i_0+1}\dashrightarrow S_{i_0+2}\dots
\] induce an infinite sequence of rational maps
\[T_1\coloneqq T\dashrightarrow T_2\dashrightarrow \dots,
\]
which are steps of a $(K_T+\Delta_T)$-MMP. Since, we are assuming termination of flips for $\mathbb Q$-factorial klt pairs in dimension $r$, we get our desired contradiction. 
\end{proof}

%
%\begin{proposition}\label{p_fdlt}
%Let $(\mathcal F,\Delta)$ be a F-dlt foliated pair on a normal projective variety $X$ such that $\mathcal F$ is algebraically integrable. Let $\pi\colon X'\to X$ be a Property $(*)$-modification and let $\mathcal F'$ be the foliation induced on $X'$. 
%
%Then every divisor in  the exceptional locus of $\pi$ is $\mathcal F'$-invariant. 
%\end{proposition}
%\begin{proof}
%
%\end{proof}

%
%\footnote{C: nothing after this footnote is needed elsewhere in the paper. I've just copied it from
%the old Section 4 (where it might have been needed if we were going to prove a basepoint free theorem).
%I've put it here just because I think it is a bit interesting, and might make a nice addition, but I'm happy
%to save it for another day.}

%!TeX root=main.tex

\section{The MMP for algebraically integrable foliations}

We begin with the following easy result: 

\begin{lemma}
\label{lem_face_containment}
Let $X$ be a normal projective variety and let  
$D$ and $A$ be $\mathbb Q$-Cartier divisors such that 
$D+tA$ is nef for all $0\leq  t \leq 1$.  

Then for all $0<t \leq 1$, we have
\[
(D+tA)^\perp \cap \overline{NE}(X) \subset (D+A)^\perp \cap \overline{NE}(X).
\]
%Moreover, if $\varphi\colon X\dashrightarrow X'$ is a $(D+A)$-trivial birational map between normal varieties and  we assume that $D'\coloneqq \varphi_*D$ and $A'\coloneqq \varphi_*A$ are $\mathbb Q$-Cartier divisors on $X'$, then 
%\[
%(D'+tA')^\perp \cap \overline{NE}(X') \subset (D'+A')^\perp \cap \overline{NE}(X').
%\]
\end{lemma}

\begin{proof}
Suppose $(D+tA)\cdot \alpha = 0$ for some $\alpha \in \overline{NE}(X)$ and $0<t\leq 1$. 
Since, by assumption, $D\cdot \alpha\ge 0$ and $(D+A)\cdot \alpha \ge 0$, 
it follows that $(D+A)\cdot \alpha=0$ and the claim follows. 
%The second claim follows immediately by the projection formula. 
\end{proof}

We now prove that termination of flips in dimension $r$ implies the  existence of a divisorial contraction or a flip associated to a $(K_{\cal F}+\Delta)$-negative extremal ray for some algebraically integrable foliation $(\cal F,\Delta)$ of rank $r$ and with log canonical singularities:

\begin{theorem}
\label{thm_existence_contraction} 
Assume termination of flips for $\mathbb Q$-factorial klt pairs in dimension $r$. 
Let $X$ be a normal projective variety, let $\mathcal F$ be an algebraically integrable foliation on $X$ of rank $r$
and let $\Delta \geq 0$ be a $\mathbb Q$-divisor such that $(\mathcal F, \Delta)$ has log canonical
%\footnote{C: I replaced F-dlt by log canonical, I think this is fine} 
singularities
and $(X, \Delta)$ is klt. Let $R\subset \overline{NE}(X)$ be a $(K_{\mathcal F}+\Delta)$-negative extremal ray.

Then there exists a contraction $c\colon X \rightarrow Y$ associated to $R$ and the following properties hold:
\begin{enumerate}
\item $Y$ is projective; and
\item if $X$ is $\mathbb Q$-factorial then $\rho(X/Y) = 1$. Moreover, if in addition 
$c\colon X \rightarrow Y$ is a flipping contraction, then the flip associated to $R$ exists.
\end{enumerate}

\end{theorem}
\begin{proof}
Let $h\colon \tilde X\to X$ be a $\mathbb Q$-factorialisation of $X$. Since $(X,\Delta)$ is klt, we may assume that $\Exc h$ has codimension at least two. 
Let $\pi\colon \overline{X} \rightarrow X$ be a Property $(*)$ modification
of $(\mathcal F, \Delta)$, which is guaranteed to exist by Theorem \ref{t_starmod}. 
We may assume that $\pi$ factors through $h$ and we denote by $\pi'\colon \overline X\to \tilde X$ the induced morphism. 
Let $\overline{\cal F}$ be the induced foliation on $\overline X$ and  let $G \ge 0$ be the associated invariant divisor. In particular, $\overline {\mathcal F}$ is induced by an equidimensional morphism $g\colon \overline X\to \overline Z$ and there exists a divisor $\Sigma_{\overline Z}$ on $\overline Z$ such that $(\overline Z,\Sigma_{\overline Z})$ is log smooth and $G$ is the pre-image, as a reduced divisor, of $\Sigma_{\overline Z}$ in $\overline X$. 
%By Proposition \ref{p_fdlt}, it follows that every divisor in  $\Exc \pi$ is $\overline{\mathcal F}$-invariant.

Let $\Gamma\coloneqq  \pi_*^{-1}\Delta$ and $\overline \Delta\coloneqq \Gamma+\sum\epsilon(E_i)E_i$, where the sum runs over all the exceptional divisors 
%$E_1,\dots,E_\ell$ 
of $\pi$. 
Thus, we may write 
\[
K_{\overline{\mathcal F}}+\overline{\Delta} = \pi^*(K_{\mathcal F}+\Delta).
\]
We may also write 
\[
K_{\overline{X}}+ \Gamma+E_0 = \pi^*(K_X+\Delta)+F_0
\]
where $E_0, F_0 \geq 0$ are $\pi$-exceptional $\mathbb Q$-divisors which do not admit any common component.
Since $\tilde X$ is $\mathbb Q$-factorial, there exists a  $\pi$-exceptional Cartier divisor $B\ge 0$ such that $-B$ is $\pi'$-ample.  Since $(X,\Delta)$ is klt, there exists  $\delta>0$ sufficiently small so 
that $(\overline{X}, \Gamma+E_0+\delta B)$ is klt.  Let $E \coloneqq E_0+\delta B$ and 
$F \coloneqq F_0+\delta B$. In particular, we have that 
\[
K_{\overline{X}}+\Gamma+E = \pi^*(K_X+\Delta)+F.
\]
Moreover, 
after possibly adding more components to $\Sigma_Z$ and $G$, 
 we may assume that $\Gamma+E \leq \overline \Delta+G$.

By the Cone Theorem \cite[Theorem 3.9]{ACSS} there exists  a nef $\mathbb Q$-Cartier divisor $H_R$ such that 
 $H_R^\perp\cap \overline{NE}(X)=R$ and such that if $D$ is a Cartier divisor  on $X$ such that 
 $R\subset D^\perp$ then 
 $(H_R+tD)^{\perp}\cap \overline{NE}(X)=R$ 
for any sufficiently small $t>0$.   
 We may write
$H_R = K_{\mathcal F}+\Delta+A$ where $A$ is an ample $\mathbb Q$-divisor. Let $\overline A\coloneqq \pi^*A$. By construction, we know that $K_{\overline{\mathcal F}}+\overline{\Delta}+\overline{A}= \pi^*H_R$ is nef
and that, for all $\lambda <1$, we have that $K_{\overline{\mathcal F}}+\overline{\Delta}+\lambda \overline{A}$ 
is not nef. We distinguish two cases. 

\medskip 

We first assume that $H_R$ is not big. Let $0\le \nu<n$ be its numerical dimension. We define $D_i\coloneqq \pi^*H_R$ for all $1\le i\le \nu+1$ and $D_i\coloneqq\overline A$ for $\nu+1<i\le n$. Then 
\[
D_1\cdot\dots\cdot D_n=(\pi^*H_R)^{\nu+1}\cdot \overline{A}^{n-\nu-1}=0
\]
and 
\[
\begin{aligned}
-(K_{\overline{\cal F}}+\overline{\Delta})\cdot D_2\cdot\dots\cdot D_n&= (\overline A - \pi^*H_R)\cdot (\pi^*H_R)^{\nu}\cdot \overline{A}^{n-\nu-1}\\
&= (\pi^*H_R)^{\nu}\cdot \overline{A}^{n-\nu}
>0.
\end{aligned}
\]
Then, by \cite[Corollary 2.28]{spicer20}, through a general point of $\overline X$, there exists a rational curve $\overline C$ which is tangent to $\overline{\cal F}$ and such that $\pi^*H_R\cdot \overline C=0$. Thus, it follows that 
\[
(K_{\overline{X}}+\overline \Delta + G+\overline A)\cdot \overline C=\pi^*H_R\cdot \overline C=0,
\] 
and, in particular, $R=\mathbb R_+[\pi_*\overline C]$. 
Since $\Gamma+E\le \overline \Delta+G$, we also have that 
\[
(K_X+\Delta)\cdot \pi_*\overline C< 0.
\] Thus, the Theorem follows from the classic base point free theorem \cite[Theorem 3.3]{KM92}. 

%\medskip 
%
%We now distinguish the two cases above. Let us first assume that $H_R$ is not big. Then there exists a curve $C_{i_0}$ as in case (i). 
%Let $\overline C$ be the strict transform of $C_{i_0}$ on $\overline X$. Since the rational map $\overline{X}\dashrightarrow \overline{X_{i_0}}$ is $(K_{\overline{\mathcal F}}+\overline{\Delta}+\overline{A})$-trivial, it follows that 
%$(K_{\overline{\mathcal F}}+\overline{\Delta}+\overline{A})\cdot \overline{C}=0$ and, therefore, $H_R\cdot \pi*\overline C=0$. 
%This implies that $R=\mathbb R_+[\pi_*C_0]$. Moreover,  we have that $(K_{\overline X}+\overline {\Delta}+G+\overline{A})\cdot \overline C=0$. 
%Since $\Gamma+E \leq \overline \Delta+G$, it follows that    

\medskip 
We now assume that $H_R$ is big. 
Since  $(\overline{\mathcal F}, \overline{\Delta})$
satisfies Property $(*)$, by Proposition \ref{p_flipstar}, we may run a  $(K_{\overline{\mathcal F}}+\overline{\Delta})$-MMP, 
with scaling of $\overline{A}\coloneqq \pi^*A$ (cf. \cite[Remark 3.10.9]{BCHM06}). In particular, 
if we define $\overline{X}_1\coloneqq \overline X$, then 
 this MMP 
with scaling produces a (possibly infinite) sequence
of $(K_{\overline{\mathcal F}}+\overline{\Delta})$-flips and divisorial contractions, 
call them $\phi_i\colon \overline{X}_i \dashrightarrow \overline{X}_{i+1}$, and a sequence
of rational numbers $1 = \lambda_1 \geq \lambda_2 \dots$ such that, if $\overline{\mathcal F}_i$ is the foliation induced on $\overline{X}_i$ and $\overline{\Delta}_i$ and $\overline{A}_i$ are respectively the strict transform of $\overline{\Delta}$ and $\overline{A}$ on $\overline{X}_i$, then $(\overline {\cal F}_i,\overline{\Delta}_i)$ satisfies Property $(*)$ and 
\[
\lambda_i\coloneqq \inf \{t\ge 0 \mid K_{\overline{\mathcal F}_i}+\overline{\Delta}_i+t\overline{A}_i \text{ is nef}\}.
\]

%one of the following holds:
%\begin{enumerate}[(i)]
%\item $H_R$ is not big and there exists a rational curve $\overline{C}_{i_0}$,  passing through the general point of $\overline{X}_{i_0}$ and which is tangent to $\overline{\cal F}_{i_0}$  and such that $(K_{\overline{\mathcal F}_{i_0}}+\overline{\Delta}_{i_0}+\overline{A}_{i_0})\cdot \overline{C}_{i_0}=0$; or
%\item $H_R$ is big and  $\lambda_{i_0}<1$.
%\end{enumerate}

%We claim that there exists a positive integer $i_0$ such that $\lambda_{i_0}<1$.
%Indeed, as in the proof of Proposition \ref{p_flipstar}, we have that each step of this MMP is also a step in a $(K_{\overline{X}}+\overline{\Delta}+G)$-MMP.
%Since we are assuming termination of flips for $\mathbb Q$-factorial klt pairs in dimension $n-1$, Special Termination
%\cite{Shokurov03} implies that 
%after finitely many flips, the flipping locus is disjoint from any log canonical centre of $(\overline{X}_i,\overline{\Delta}_i+G_i)$, where $G_i$ is the strict transform of  $G$ on $\overline X_i$.  
%In particular, since $\Gamma+E \leq \overline \Delta+G$ and the difference is supported in the union of log canonical centres of $(\overline X, \overline \Delta+G)$, it follows that, for any sufficiently large $i$, the map $\phi_i$
%is also a step in a $(K_{\overline{X}}+\overline{\Gamma}+E)$-MMP. Since $(\overline X,\overline{\Gamma}+E)$ is klt and  $\overline A$ is big,  \cite[Corollary 1.4.2]{BCHM06} implies the claim. 

By Lemma \ref{l_termination*}, there exists a positive integer $i_0$ such that $\lambda_{i_0}<1$.
After possibly replacing $i_0$ by a smaller number, we  may assume that $i_0$ is the smallest positive integer such that $\lambda_{i_0}<1$.  
Set $\lambda \coloneqq \lambda_{i_0}$, 
$\overline{X}'\coloneqq \overline{X}_{i_0}$
and let $\phi\colon \overline{X} \dashrightarrow \overline{X}'$ be the induced birational map.
Let
$\overline{\mathcal F}' \coloneqq \phi_*\overline{\mathcal F}$, $\Gamma'\coloneqq \phi_*\Gamma$, $\overline{\Delta}' \coloneqq \phi_*\overline{\Delta}$,
$E' \coloneqq \phi_*E$ and 
$\overline{A}' \coloneqq \phi_*\overline{A}$. By the negativity lemma (cf. \cite[Lemma 3.39]{KM92}), it follows that $(\overline{X}', \Gamma'+E')$
is klt.

By definition of the MMP with scaling and by our choice of $i_0$, 
 we have that $K_{\overline{\mathcal F}'}+\overline{\Delta}'+t\overline{A}'$
is nef for all $\lambda \leq t \leq 1$ and that  each step of this MMP  up until $\overline{X}'$ 
is $(K_{\overline{\mathcal F}}+\overline{\Delta}+t\overline{A})$-negative and  $(K_{\overline{\cal F}}+\overline {\Delta} +\overline{A})$-trivial. Thus, 
$\phi_*\pi^*H_R = K_{\overline{\mathcal F}'}+\overline{\Delta}'+\overline{A}'$ is nef.
%Since $\phi_*\pi^*H_R$ is a limit of nef divisors, it follows from 
By Lemma \ref{lem_face_containment}, we have a containment
\begin{align}
\label{face_containment}
(K_{\overline{\mathcal F}'}+\overline{\Delta}'+t\overline{A}')^\perp\cap \overline{NE}(\overline{X}') \subset 
(\phi_*\pi^*H_R)^\perp\cap \overline{NE}(\overline{X}')
\end{align}
for all $\lambda < t \leq 1$.

Fix a rational number $\lambda'$ such that $\lambda< \lambda'<1$ and  a sufficiently small rational number $s>0$ such that 

\begin{itemize}
\item $K_{\cal F}+\Delta+\lambda' A$ is big;
\item if we set $A_0\coloneqq (1-\lambda' )A - s(K_X+\Delta)$, then $A_0$ is ample, the stable base locus of
$H_R-A_0$ coincides with $\mathbb B_+(H_R)$ 
and $H_R-A_0$ is positive on every extremal ray of  $\overline{NE}(X)$ except $R$;

\item $s < \frac{1}{2m\dim X}$ where $m$ is the Cartier index of 
$K_{\overline{\mathcal F}'}+\overline{\Delta}'+\lambda'\overline{A}'$; and

\item if we set $K\coloneqq (K_{\overline{\mathcal F}}+\overline{\Delta}+\lambda'\overline A)+s
(K_{\overline{X}}+\Gamma+E )$, then $K$ is big and  $\phi$ is $K$-negative.
\end{itemize}

Set $K' \coloneqq \phi_*K$.  %As in the proof of Proposition \ref{prop_adjoint_*_MMP}, 
By \cite{Kawamata91}, every $(K_{\overline{X}'}+\Gamma'+E')$-negative extremal ray is spanned 
by a curve $C$ such that
$-(K_{\overline{X}'}+\Gamma'+E')\cdot C \leq 2\dim X$. Thus, since 
$K_{\overline{\mathcal F}'}+\overline{\Delta'}+\lambda'\overline A'$ is nef and big
by our choice of $\lambda'$, we may run a $K'$-MMP which is $(K_{\overline{X}'}+\Gamma'+E')$-negative
and $(K_{\overline{\mathcal F}'}+\overline{\Delta}'+\lambda' \overline{A}')$-trivial. Call this MMP  
$\psi\colon \overline{X}' \dashrightarrow \overline{X}''$ and 
let
$\overline{\mathcal F}'' \coloneqq \psi_*\overline{\mathcal F}'$, 
$\Gamma''=\psi_*\Gamma'$, 
$\overline{\Delta}'' \coloneqq \psi_*\overline{\Delta}'$,
$E'' \coloneqq \psi_*E'$,  
$\overline A'' = \psi_*\overline A'$ and $K''\coloneqq \psi_* K'$. 

We have that $(\overline{X}'', \Gamma''+E'')$ is klt and 
$\frac 1 sK'' - (K_{\overline{X}''}+\Gamma''+E'')$ is big and nef. 
Thus, by the classic base point free theorem \cite[Theorem 3.3]{KM92}, it follows that $K''$ is semi-ample.

%Let $c\colon \overline{X}'' \rightarrow Y$ be the morphism associated to $K''$.  Note
%that the Cartier index of $K_{\overline{\mathcal F}''}+\overline{\Delta}''+(1-t)A''$
%remains $m$, and so (arguing as in the proof of 
%Proposition \ref{prop_adjoint_*_MMP}) we see that $c$ only contracts curves
%which are both $K_{\overline{X}''}+\overline{\Delta}''+E''$-trivial and 
%$K_{\overline{\mathcal F}''}+\overline{\Delta}''+(1-t)A''$-trivial.
%We remark that since $(\overline{X}'', \overline{\Delta}''+E'')$ is klt that
%$c$ may also be realised as a $K_{\overline{X}''}+\Theta'$-negative contraction
%for some klt pair $(\overline{X''}, \Theta'')$.  In particular, if $M$ is $c$-numerically trivial,
%then $M = c^*N$ for some $\mathbb Q$-Cartier divisor $N$ on $Y$.

Recall that, by our choice of  $\lambda'$ and $s$,  the $\mathbb Q$-divisor 
$A_0=(1-\lambda' )A - s(K_X+\Delta)$ is ample and  the stable base locus of
$H_R-A_0$ coincides with $\mathbb B_+(H_R)$.
We may write $K = \pi^*(H_R-A_0)+sF$ and by Nakamaye's theorem (cf. \S \ref{s_preliminaries}) it follows that 
the restricted base locus of $K$ is exactly the union of $\Supp F$ with $Z$, the preimage
of the augmented base locus of $H_R$.  
Thus, Lemma \ref{l_exceptional} implies that the divisorial part of $\Supp F\cup Z$ is contracted by $\psi \circ \phi$, and since
$\Exc \pi' = \Supp F$ and $h$ is small, it follows that the induced map $f\colon X \dashrightarrow \overline{X}''$ 
is a birational contraction. In particular, we have that 
$\Gamma''=\overline{\Delta}''=f_*\Delta$. 
Since $\tilde X$ and $\overline{X}''$ are both $\mathbb Q$-factorial, it also follows that if the augmented base locus of $H_R$ does not 
contain any divisor then
$\rho(\tilde X) = \rho(\overline{X}'')$ and, otherwise, $\rho(\tilde X)-1 = \rho(\overline{X}'')$.
Set $X^+\coloneqq \overline{X}''$.  

\medskip

Next, we claim that $\varphi\coloneqq \psi\circ \phi\colon \overline X\dashrightarrow {\overline X}''$ is $(\pi^*H_R)$-trivial.  As noted earlier, this holds for $\phi$
by definition of the MMP with scaling.  It remains to prove that $\psi$ is $(\varphi_*\pi^*H_R)$-trivial. 
Let $\psi_i\colon \overline{X}'_i \dashrightarrow \overline{X}'_{i+1}$ 
be the steps of the MMP composing $\psi$ for $i=0,\dots,k-1$, 
where $\overline{X}'_0\coloneqq \overline{X}'$ and $\overline{X}'_k \coloneqq \overline{X}''$.
Let $\overline{\mathcal F}'_i$ be the induced foliation on $\overline{X}'_i$ and 
let $\overline{\Delta}_i'$, $\overline{A}_i'$ and $(\phi_*\pi^*H_R)_i$
be the strict transform of  
$\overline{\Delta}$, $\overline{A}$ and $\pi^*H_R$  
on $\overline X'_i$, respectively. 
By construction,
$\psi_i$ is $(K_{\overline{\mathcal F}'_i}+\overline{\Delta}'_i+\lambda' \overline A_i')$-trivial.  
We claim that we have a containment for all $\lambda <t \leq 1$ and $i$.
\begin{align}
\label{face_containment_2}
(K_{\overline{\mathcal F}_i'}+\overline{\Delta}_i'+t\overline{A}_i')^\perp \cap \overline{NE}(\overline{X}_i')
\subset (\phi_*\pi^*H_R)_i^\perp \cap \overline{NE}(\overline{X}_i').
\end{align}
Supposing this containment we may conclude immediately.  To prove the claim we argue by induction on $i$.   
When $i = 0$ this is
exactly 
the containment \eqref{face_containment}.  
Supposing the claim for $i\ge 0$, it follows that $\psi_i$ is $\overline A_i'$-trivial. In particular, both 
$(K_{\overline{\mathcal F}_{i+1}'}+\overline{\Delta}_{i+1}'+\lambda\overline{A}_{i+1}')$
and
$(K_{\overline{\mathcal F}_{i+1}'}+\overline{\Delta}_{i+1}'+\overline{A}_{i+1}')$
are nef and Lemma \ref{lem_face_containment} implies \eqref{face_containment_2}. Thus, our claim follows.

\medskip

%By perhaps perturbing $A_0$ slightly we may assume that 
%$f_*(H_R-A_0)$ is ample. Indeed, let $H$ be an ample divisor on $\overline{X}''$.
%Then for all $\epsilon>0$ we have $f_*(H_R-A_0)+\epsilon H$ is ample.
%For $\epsilon>0$ sufficiently small we have that $A_0-\epsilon f_*^{-1}H$ is ample
%and $\psi\circ \phi$ is still an MMP for  $\pi^*(H_R-A_0+\epsilon f_*^{-1}H)+sF$.
%So by replacing $A_0$ by $A_0-\epsilon f_*^{-1}H$ we may assume that $f_*(H_R-A_0)$
%is ample.

We now show that there exists a contraction $c_R\colon X \rightarrow Y$ associated to $R$.
  If  $(K_X+\Delta)\cdot C\le 0$ for any $C\in R$, then the claim follows from the   
  base point free theorem (cf. \cite[Theorem 3.3]{KM98}). 
Thus, we may assume that  $R$ is $-(K_X+\Delta)$-negative and, therefore,  there exists  
$c>0$ such that $c(K_X+\Delta)\cdot \xi=A_0\cdot \xi$ for all $\xi\in R$.  
By our choice of $H_R$, we have that  if $m$ is a sufficiently large positve integer and  
\[
\tilde{H}_R \coloneqq (c(K_X+\Delta)-A_0)+mH_R
\]
then $\tilde H_R$ is nef and $\tilde H_R\cap \overline{NE}(X)=R$.   In particular, we have that
\[
\pi^*H_R^\perp\cap {\overline {NE}}(\overline X)=\pi^*\tilde{H}_R^\perp\cap {\overline {NE}}(\overline X). 
\]
and since $\varphi$ is $(\pi^*H_R)$-trivial, it follows that it is also $(\pi^*\tilde{H}_R)$-trivial. 
In particular, $f_*\tilde H_R$ is nef. 
Since $f_*(K_X+\Delta) = K_{X^+}+\Delta^+$, where $\Delta^+\coloneqq \overline{\Delta}''$,  
we have
that 
\[
f_*\tilde H_R - c(K_{X^+}+\Delta^+) = f_*(H_R-A_0)+(m-1)f_*H_R
\]
is big and nef, since $f_*(H_R-A_0) =\varphi_*K$ is nef and $f_*H_R$ is big and nef.
Since $(X^+,\Delta^+)$ is klt and $\mathbb Q$-factorial, we may therefore apply the 
base point free theorem to conclude 
that $f_*\tilde{H}_R$ is semi-ample.  It follows that $\tilde{H}_R$ itself is 
semi-ample and, therefore, the morphism associated to $\tilde H_R$ is the contraction 
$c_R\colon X \rightarrow Y$ associated to $R$ and our claim follows. 
Note that, by construction, $Y$ is projective  and $\rho(\tilde X/Y)=1$. 

\medskip

If $X$ is $\mathbb Q$-factorial, it follows that $\rho(X/Y)=1$. 
We now assume, in addition, that $c_R\colon X\to Y$ is a flipping contraction. 
Then $\rho(X)=\rho(X^+)$ and the morphism associated to $f_*\tilde H_R$ 
defines a birational morphism $c^+\colon X^+\to Y$. Thus, $\rho(X^+/Y)=1$. 
Since $(H_R-A_0)\cdot \xi <0$ for any curve $\xi$  contracted by $c$ 
and  
\[
f_*(H_R-A_0)\cdot \xi^+=\phi_*K\cdot \xi^+\ge 0
\] for any curve $\xi^+$ contracted by $c^+$, 
it follows that $f\colon X\dashrightarrow X^+$ is the flip associated to $R$.  
\end{proof}

\begin{remark}\label{r_mmp}
Using the same notation as in the proof of the Theorem, assume that $c_R$ is a flipping contraction. Then the map $\psi\colon \overline X'\dashrightarrow \overline X''=X^+$ is $(K_{\overline{\mathcal F'}}+\overline{\Delta}'+\lambda'\overline{A}')$-trivial and  $K_{\overline{\mathcal F}''}+\overline{\Delta}''+\lambda'\overline{A}''$ is ample over $Y$. Thus, Lemma \ref{l_amplemodel} implies that $\psi$ is a $(K_{\overline{\mathcal F}'}+\overline{\Delta}')$-trivial morphism. Note also that $\overline{X}'$ is klt and $\mathbb Q$-factorial, the induced morphism $\overline {X}'\to Z$ is equidimensional, $(\overline{\mathcal F}',\overline{\Delta}')$ satisfies Property $(*)$ and
every $\overline{\cal F}'$-invariant divisor in the exceptional locus of $\psi$ is contained in the support of the divisor associated to 
$(\overline{\mathcal F}',\overline{\Delta}')$.

Thus, if $f\colon X\dashrightarrow X^+$ is a $(K_{\cal F}+\Delta)$-flip,  $\pi\colon \overline X\to X$ is a Property $(*)$-modification for $(\cal F,\Delta)$ and
$\cal F^+$ is the induced foliation on $X^+$, then there exists  a Property $(*)$-modification  $\pi^+\colon \overline X^+\to X^+$ for $(\cal F^+,f_*\Delta)$, such that if $\overline {\cal F}$ is the foliation induced on $\overline X$ and we write
\[
K_{\overline{\mathcal F}}+\overline{\Delta} = \pi^*(K_{\mathcal F}+\Delta),
\]
then the induced map $\overline f\colon \overline X\dashrightarrow \overline X^+$ is a sequence of steps of a $(K_{\overline{\cal F}}+\overline{\Delta})$-MMP and if $\overline {\cal F}^+$ is the induced foliation on $\overline X^+$ then 
\[
K_{\overline{\mathcal F}^+}+f^+_*\overline{\Delta} = \pi^{+*}(K_{\mathcal F^+}+f_*\Delta).
\]
\end{remark}

%
%\iffalse
%It follows that if $Z$ contains no divisors then
%$f\colon X \dashrightarrow X^+$
%is the $K_{\mathcal F}+\Delta$-flip associated to $R$, and if $Z$ contains a divisor then $f$
%is in fact a morphism and gives the divisorial contraction associated to $R$. 
%
%
%
%Since each step of our MMPs is $\pi^*H_R$ trivial, it follows that 
%$f_*H_R$ is nef.  We claim that $f_*H_R$ is zero on exactly one extremal ray, 
%call it $R'$.
%Supposing the claim we see that $H_R$ (equivalently $f_*H_R$) 
%is semi-ample as follows.  If $R$ is $K_X+\Delta$-negative,
%then $H_R-\epsilon (K_X+\Delta)$ is ample for $0<\epsilon \ll 1$ and we may apply 
%the base point free theorem, \cite[Theorem 3.3]{KM98}, to conclude.
%Otherwise $R'$ is $K_Y+f_*\Delta$-negative.  As observed above $(Y, f_*\Delta = c_*\overline{\Delta}''+E'')$
%is klt and so again $f_*H_R - \epsilon (K_Y+f_*\Delta)$ is ample for $0<\epsilon \ll 1$, and so 
%by the base point free theorem $f_*H_R$ is semi-ample.
%\fi
%
%\medskip
%
%We now assume that $H_R$ is not big, it is not hard to deduce that in fact
%$R$ is a $K_X+\Delta$-negative extremal ray\footnote{P: we should add a proof}, in which case we may conclude by
%applying the classical version of the contraction theorem.

\medskip

We are now ready to prove our main Theorem: 

\begin{proof}[Proof of Theorem \ref{t_main}]
By Theorem \ref{thm_existence_contraction}, all the divisorial contractions and flips required to run an
MMP exist. It remains to show that the MMP terminates. To this end, we show that there is no infinite sequence of flips.

%We first assume that $(\mathcal F, \Delta)$ admits Property $(*)$. Then the foliation $\cal F$ is induced by a morphism 
%$X \rightarrow Z$. Let $G$ be the associated divisor. 
%Then, as in the proof of Proposition \ref{p_flipstar}, it follows that  each step in the MMP is a step in a $(K_X+\Delta+G)$-MMP over $Z$.
%Since $\dim (X/Z) \leq n-1$, there is no infinite sequence of $(K_X+\Delta+G)$-flips by assumption, and so the 
%$(K_{\mathcal F}+\Delta)$-MMP must terminate.
%We now consider the general case. 
Suppose there exists a sequence of $(K_{\cal F}+\Delta)$-flips
\[
X_1\coloneqq X\dashrightarrow X_2\dashrightarrow \dots
\]
Let 
%$Y_i$ be the base of the flip $X_i \dashrightarrow X_{i+1}$ and let
$\cal F_i$ be the induced foliation on $X_i$ and let $\Delta_i$ be the strict transform of $\Delta$ on $X_i$. 
Then, by Remark \ref{r_mmp}, there exist Property $(*)$ modifications $\pi_i\colon \overline {X}_i\to X_i$ of $(\cal F_i,\Delta_i)$ such that, if
we write
\[
K_{\overline{\mathcal F}_1}+\overline{\Delta}_1 = \pi_1^*(K_{\mathcal F}+\Delta),
\]
then
 the induced maps
\[
\overline{X}_1\dashrightarrow \overline{X}_2\dashrightarrow \dots
\]
are a sequence of steps of a $(K_{\overline{\mathcal F}_1}+\overline{\Delta}_1)$-MMP and, therefore, by Lemma \ref{l_termination*}, it terminates. 
Thus, also the sequence of $(K_{\cal F}+\Delta)$-flips terminates, as claimed. 
%Let $\mu\colon X'_1 \rightarrow X_1$ be a Property $(*)$ modification of $(\cal F,\Delta)$,  which is guaranteed to exist by Theorem \ref{t_starmod}.
%Let $\psi_1 \colon X'_1 \dashrightarrow X'_2$ be a $(K_{\mathcal F'_1}+\Delta'_1)$-MMP over $Y_1$.
%%By Proposition  Proposition \ref{p_flipstar},
%Note that, since  $(\mathcal F'_1, \Delta'_1)$
%admits Property $(*)$, by the previous case, it follows that this MMP terminates. 
%Consider the birational contraction $c\colon X'_2 \dashrightarrow X_2$.  
%By assumption, $K_{\mathcal F_2}+\Delta_2$ is ample over $Y_1$, and so $X_2$ is the ample model
%of $K_{\mathcal F'_2}+\Delta'_2$ over $Y_1$, in particular, Lemma \ref{l_amplemodel} implies that  $c$ is a morphism.
%
%Continuing inductive we produce a sequence of MMPs $\psi_i\colon X'_i \dashrightarrow X'_{i+1}$, 
%which all form steps of  $K_{\mathcal F'_1}+\Delta'_1$-MMP. However, since $(\mathcal F'_1, \Delta'_1)$
%has Property $(*)$ this process must eventually terminate, and therefore the sequence of flips
%$X_i \dashrightarrow X_{i+1}$ must also terminate.
\end{proof}

\section{Some open problems}
We conclude with some open problems related to the main results in this paper. 
First, we expect the following generalisation of the base point free theorem to be true:
\begin{conjecture}
Let $X$ be a normal projective variety, let $\mathcal F$ be an algebraically integrable foliation on $X$
and let $\Delta \geq 0$ be a $\mathbb Q$-divisor so that $(\mathcal F, \Delta)$ has F-dlt singularities
and $(X, \Delta)$ is klt. Let $A$ be an ample $\mathbb Q$-divisor on $X$ such that $L=K_{\cal F}+\Delta+A$ is a nef $\mathbb Q$-divisor. 

Then $L$ is semi-ample. 
\end{conjecture}

In addition, we anticipate that, as with foliations on a three-dimensional variety (e.g. see \cite[Theorem 1.6 and Theorem 11.3]{CS18}), the Minimal Model Program will yield some interesting implications on  foliated singularities:
\begin{conjecture}
Let $X$ be a $\mathbb Q$-factorial projective variety, let $\mathcal F$ be an algebraically integrable foliation on $X$
and let $\Delta \geq 0$ be a $\mathbb Q$-divisor so that $(X,\Delta)$ is klt 
and one of the following holds:
\begin{enumerate}
\item $(\cal F,\Delta)$ is F-dlt, or 
\item $(\cal F,\Delta)$ is canonical. 

\end{enumerate}
Then there exists a morphism $f\colon X \rightarrow Y$ which induces 
$\cal F$. %has non-dicritical singularities. 
\end{conjecture}

Finally, we expect the following: 

\begin{conjecture}
Let $X$ be a normal projective variety, let $\mathcal F$ be a foliation on $X$ 
and let $\Delta \geq 0$ be a $\mathbb Q$-divisor on $X$ so that $(X,\Delta)$ is klt, $(\cal F,\Delta)$ admits log canonical singularities and $K_{\cal F}+\Delta$ is not pseudo-effective. 

Then we may run a $(K_{\cal F}+\Delta)$-MMP. 
\end{conjecture}
Note that we are not assuming that $\cal F$ is algebraically integrable. On the other hand, \cite{CP19} implies  that if $\cal F$ is a foliation such that $K_{\cal F}$ is not pseudo-effective then the algebraic part of $\mathcal F$ is non-trivial, i.e. there exist a dominant map $\sigma\colon X\dashrightarrow Y$ and a foliation $\cal G$ on $Y$ such that $\dim Y<\dim X$ and $\cal F=\sigma^{-1}\cal G$.

\bibliography{math}
\bibliographystyle{alpha}

\end{document}